\input amstex
\documentstyle{amsppt}


\magnification=1100


\vsize=8.3in


\baselineskip 16pt


\parindent 20 pt
\NoBlackBoxes

\TagsOnRight

\define \a{\alpha}
\define \be{\beta}

\define \g{\gamma}
\define \G{\Gamma}

\define \vp{\varphi}
\define \Vp{\Phi}

\def\cupl{\operatornamewithlimits{\bigcup}\limits}

\define \ov{\overline}

\define \ep{\endproclaim}

\define \1{^{-1}}
\define \2{^{-2}}

\centerline{\bf Counting rational maps onto surfaces}
\centerline{\bf and fundamental groups}

\bigskip

\centerline{\bf T. Bandman
\footnote{Supported by the Ministry
 of Absorption (Israel), the
Israeli Science Foundation (Israeli Academy of
Sciences, Center of Excellence Program), the Minerva
Foundation (Emmy Noether Research Institute of
Mathematics) }}
\centerline {Department of Mathematics}
\centerline {Bar Ilan University}
\centerline {Ramat-Gan, Israel}
\centerline{bandman\@macs.biu.ac.il}
\centerline{\bf A.Libgober \footnote{Supported by NSF grant}}
\centerline {Department of Mathematics}
\centerline {University of Illinois at Chicago}
\centerline {851 S.Morgan, Chicago, Ill, 60607}
\centerline{libgober\@math.uic.edu}

\topmatter
 
 \leftheadtext{T. Bandman and A.Libgober}
 \rightheadtext{ Counting rational maps onto surfaces }


\abstract
 We consider the class of quasiprojective varieties admitting a 
 dominant morphism onto a curve with negative Euler 
characteristic. The existence of such a morphism
is a property of the fundamental group.
We show that for a variety in this class the number of 
maps onto a hyperbolic curve or surfaces 
can be estimated in terms of the numerical invariants 
of the fundamental group. We use this estimates 
to find the number of biholomorphic automorphisms of 
complements to some arrangements of lines.
\endabstract
\endtopmatter

\subheading {1.Introduction}
 In this paper we study quasiprojective varieties $X$
admitting a dominant morphism onto a curve with negative Euler 
characteristic. A fundamental result of D.Arapura (\cite{1})
shows that such maps can be described in terms of the moduli spaces
of rank one local systems with non vanishing cohomology.
His work was built on previous extensive studies going back 
to DeFranchis and, more recently, to A.Beauville (\cite{2}), 
Catanese, Green-Lazarsfeld and Deligne-Simpson. 

The moduli spaces of local systems with non-vanishing $H^1$ 
can be described in terms 
of the invariants of the fundamental group of $X$ alone.
These invariants, in fact depending only on 
the quotient of the fundamental group by its 
second commutator, are the algebraic subvarieties of 
$Hom(\pi_1(X),{\Bbb C}^*)$.
They were considered 
also in (\cite{3} in the case when $X$ is 
a complement to a plane curve and called the characteristic varieties.
A consequence of the result of Arapura is that 
the number of possible dominant maps
with connected fibers
from a quasiprojective variety $X$ onto hyperbolic 
(i.e. having negative Euler characteristic) curves
is equal to the number $n(X)$
 of those irreducible components of positive dimension in
 the characteristic variety of $\pi_1(X)$ which 
contain the identity of $Hom(\pi_1(X),{\Bbb C}^*)$ 
(the number of such $i-$dimensional
 components we denote by $n_i(X)$).

 Here we use the characteristic varieties 
 to find effective bounds for the number
 of maps and targets which are surfaces of hyperbolic type.
For two quasiprojective varieties of hyperbolic type the number $R(X,Y)$
 of dominant morphisms of $X$ into $Y$ is finite (\cite{4}, \cite{5},
\cite{6},
\cite{7}). The number of these maps and the number of
targets (\cite{8}) was a subject of numerous works. 
In the case of automorphisms, for a projective variety
$X$ with at most terminal singularities, nef canonical bundle $K(X)$ 
and having an index $r$ and the dimension $n$,
 the number of birational automorphisms
has polynomial bound  of the form  $f(n)(rK(X)^n)^{g(n)}$ (\cite{9}).
For a minimal surface  one has:
$\# Aut(X)\leq CK(X)^2$ (\cite{10},\cite{11}); for a smooth
 hyperbolic curve $X$ of genus $g$ with $k$ punctures
$\# Aut(X)\leq A(g,k)$, where 
 $$ A(g,r)\le \left \{ \matrix\format \l&&\quad\l\\
5k  &&\text {if} \   g=0\\
6(2g+k-2) && \text {if} \ k>0, g\geq1\\
84(g-1) && \text {if}  \  k=0 \endmatrix\right.\tag 1$$
(see \cite{12}).

In the case when $dim X=1$, 
the number of maps of a compact curve of genus $g$
 onto other hyperbolic curves (up to isomorphism of the image) is bounded by

$$T_c(g)=(g-1)2^{2g^2-2}(2^{2g^2-1}-1)\tag 2$$
\noindent (\cite{13}, \cite{14},\cite{15}). 
This bound cannot be significally improved (\cite{14}).

Non-effective bounds for $R(X,Y)$ were obtained in 
\cite{16}, \cite{17} in 
the case when $X,Y$ are projective varieties of general type
having at most canonical singularities
and $K(X),K(Y)$ are nef.
The number of targets for threefolds was found in (\cite{18})
and the bounds for some special classes of pairs of $(X,Y)$ 
were obtained  in (\cite{19},\cite{20}.
The case of canonically polarized varieties was considered in (\cite{21}).

An effective estimate for the number of holomorphic maps 
of a quasiprojective $X$ with the target 
being an arbitrary hyperbolic curve is 
available in the case 
when $dim X=1$. If $X$ is smooth, has genus $g$ and $r$ punctures
then the number of its dominant maps with hyperbolic targets
 is at most
 $$T(g,r)=T_c(g)+T_1(g,r)+T_0(g,r),\tag 3$$
where $T_c(g)$ is given by (2),
$$ T_1(g,r)= \left \{ \matrix\format \l&&\quad\l\\
r2^{3g^2-1}(2g+r)^{2g^2-2}[\zeta(2g^2-2)\zeta (2g^2-3)]  &&
\text {if} \   g>1\\
 r^2ln(r) && \text {if} \  g=1\\
0 && \text {if}  \  g=0 \endmatrix\right.\tag 4$$
and 
$$T_0(g,r)=[3(2g+r-1)]^{(2g+r)}\tag 5$$
Here $T_1(g,r)$ (resp. $T_0(g,r)$)
 is the bound for the maps onto punctured tori (\cite{14})
(resp. the number of possible maps onto punctured plane (\cite{22}));
$\zeta$  is the Riemann $\zeta$-function).

Previous estimates in the above problems about the 
maps of $X$ with hyperbolic targets depend on the numerical  
data of $X$.
The purpose of this paper is to describe the estimates 
in such kind of problems rather 
in terms of the {\it fundamental group} 
and find explicit bounds for the number of maps
onto the curve and surfaces for which the characteristic varieties
have components of  positive dimension.

Sections 2 contain preliminaries.
In  section 3
 we consider maps of projective varieties (Theorems 1 and 2) and
in section 4 we
find the bounds for
the number of maps between
 smooth affine surfaces $S$
 with $n(S)>0$ (Theorems 3,4,5). Here is a sample of 
our results:
\proclaim{Theorem 4}
   Let $S$ be an affine surface of general type with $n(S)\ge 0,$
 Then
$$\#Aut(S)\le\min_{\phi,\nu _{\G}\ne 0}\{\nu _{\G}\cdot
L(\Gamma)\cdot L_{_{\phi}}\}.$$
where   $\nu _{\G}$ is the length of an orbit of a component 
of characteristic variety corresponding to a map $\phi: 
S \rightarrow {\G}$, \  $L(\Gamma)=\# Aut(\G),$
 and $F_{\phi}$ is a generic fiber of  $\phi$ 
and $L_{\phi}= \max\{\# Aut (F_\phi)\}$ for all such $\phi.$ 
In particular 
$$\#Aut(S)\le \min_{n_i(S)\ne 0} \{252 n_i(S)\cdot e(S)\}.$$
\ep

As an application of this theorem, in the appendix we calculate the 
order of the {\it biholomorphic} automorphisms
to the complements to line arrangements.
In particular the order of the automorphism of  
the complement to the Ceva arrangement
is equal to $120$ and hence
contains non linear automorphism (cf. example 1). On the other hand all 
biholomorphic maps of the complements to Hesse arrangement 
are linear. 

 We have to acknowledge that for projective varieties the
estimates are so big, that their  main importance is just the information
about the
invariants of the variety, which are responsible  for the existence and the
number
of these maps. It seems that $n(X)$ is just one of these invariants.

 For affine varieties  (section 4) the situation is different.
 Theorem 3 provides not only 
the effective estimate of the number of maps into curves, but
 shows that the main
 invariant
in this game is the fundamental group. This theorem is an answer
to a long - standing
question, formulated originally in terms of the holomorphic functions
(see Remarks after Theorem 3). The estimate  for the number of automorphisms
for affine surfaces (Theorem 4) could be obtained without $n(S),$ but then it
 would
 include the  additional factor of type
   $T(g,r),$ which is exponential, instead
 of polynomial one.

 Our estimates include the function $T(g,r)$
though we never use its explicit form.
Roughly speaking, in projective case we prove that if $T$ is the best
 estimate for the number of maps into hyperbolic curves, then $T^2$
(with some modifications) is the estimate for the number of dominant
maps into surfaces with special fundamental groups.

The main idea used in this paper is the following.
If  $S_1 \rightarrow S_2$ is a map between the surfaces 
and $S_2,$ is equipped with a fibration,
 then $S_1$
is equipped by a fibration as well. Any map of a  fibration with base
$\Gamma _1$  and a general fiber $F_1$  into a fibration with base
$\Gamma _2$ and a general fiber $F_2$ provides the maps of $\Gamma _1$
into $\Gamma _2$ and $F_1$ into $F_2.$ Since all these curves are
hyperbolic, there exist estimates for the  number of maps.
The number of fibrations is given  in terms of the fundamental group.

Finally, we use the following terminology: a variety $X$ 
is of general (hyperbolic) type
if it has a smooth compactification $\ov X$ such that :

a) $X=\ov X-D,$ where $D$ is a normal crossing divisor;

b) linear system $|m(K(\ov X)+D)|$ provides a birational map for some $m.$

\bigskip

\subheading{2. Preliminaries: characteristic variety of the fundamental group
and maps onto the curves}

\bigskip
In this section we shall describe the characteristic varieties 
of the fundamental
group in terms of which we estimate of the order of the
automorphisms groups. They depends only on the 
quotient $\pi_1'(X)/\pi_1''(X)$ of the fundamental group by its second 
commutator. For additional details cf. (\cite{3},\cite{23}).

Let $G$ be a finitely generated, finitely presented group such
that the abelianization $G/G' \ne 1$. Let $r$ be a number of
generators of $G/G'$ and let $s$ be a surjection:
${\Bbb Z}^r \rightarrow G/G'$. We shall consider the exact
sequence:

$$0 \rightarrow G'/G'' \rightarrow
 G/G'' \rightarrow G/G' \rightarrow 0$$

\noindent where $G''=[G',G']$ is the second commutator. Since 
the left group in this sequence is abelian, a lifting the elements
in $G/G'$ to elements in $G/G''$ yields the action of $G/G'$ on
$G'/G''$ and the module $M_G=G'/G'' \otimes {\Bbb C}$ over the
group ring ${\Bbb C}[G/G']$ of $G/G'$. The characteristic variety $C_G$ of
$G$ is the support of this module. This is the subset of
$Spec {\Bbb C}[G/G']$ consisting of prime ideals in
$\wp \subset {\Bbb C}[G/G']$ such that $M_G \otimes_{{\Bbb C}[G/G']}
 ({\Bbb C}[G/G']/\wp) \ne 0$.

The surjection $s$ allows to view $Spec {\Bbb C}[G/G']$ as a subset
in the torus ${{\Bbb C}^*}^r$. $C_G$ has canonical filtration defined
as follows. Let $\Phi: {\Bbb C}[G/G']^m \rightarrow
{\Bbb C}[G/G']^n \rightarrow M_G \rightarrow 0$ be a presentation
of the module $M_G$ with $n$ generators and $m$ relations.
Let $F_i(G)$ be the ideal in ${\Bbb C}[G/G']$ generated by
$(n-i+1) \times (n-i+1)$ minors of the matrix of $\Phi$
(Fitting ideal of the module $M_G$).
Let $C_G^k$ be the (reduced) zero set of $F_k(G)$.
We have the inclusions:
$$... \subseteq C^{k+1}_G \subseteq C^k_G \subseteq ... \subseteq C^1_G
\subseteq Spec {\Bbb C}[G/G']$$

\noindent All affine
varieties $C_G^k$ and ideals $F_k(G)$ are invariants of the
fundamental group depending only on $G/G''$.

Let $X$ be a connected CW-complex. A local
system of rank one on $X$ is a homomorphism $\pi_1(X) \rightarrow
{\Bbb C}^*$ i.e. a character of the fundamental group $G=\pi_1(X)$. 
Non trivial
rank one local systems exist if and only if $G/G' \ne 0$. Clearly
the group of characters of $G$ can be identified with $Spec {\Bbb
C}[G/G']$. The above filtration on $Spec {\Bbb C}[G/G']$ 
has a description in terms of local systems 
as follows.

The cohomology $H^i(X,\rho)$ of a local system $\rho$ can be defined as
the cohomology of a chain complex:
$$... \rightarrow C_*(\tilde X,\Bbb C) \otimes_{\pi_1(X)} {\Bbb C}_{\rho}
\rightarrow ...$$
where $\tilde X$ is the universal cover of $X$, $C_*(\tilde X,{\Bbb C})$
is a $\Bbb C$-vector space of chains on $\tilde X$ 
considered as a $\pi_1(X)$-module with $\pi_1(X)$-action
coming from the action on $\tilde X$ and ${\Bbb C}_{\rho}$ is $\Bbb C$
with the structure of ${\Bbb C}[G/G']$-module via $\rho$.
The set $C_G^k$ coincides with the set of characters $\rho$
with the property $dim H^1(X,\rho) \ge k$.

In the case when $X$ is a quasi-projective variety and $G=\pi_1(X)$,
the varieties $C_G^k$ have a
remarkably simple structure 
(cf. \cite{1}): each $C^k_G$ is a union of translated subgroups
of $ Spec {\Bbb C}[G/G']$. Moreover, one has the following.
Let us call a map $f: X \rightarrow C$ {\it admissible} if 
there is an extension $\bar f: \bar X \rightarrow \bar C$ 
with smooth $\bar X, \bar C$  such that $\bar f$ has 
connected fibers (cf. \cite{1}). 
For any quasi-projective $X$, there exist unitary characters
$\rho_j^{\prime}$, the torsion characters $\rho_i$, the curves
$C_i$ and the admissible maps $f_i: X \rightarrow C_i$ such that
$$C^1_{\pi_1(X)}=\bigcup \rho_i f^*_i(H^1(C_i)) \cup \bigcup \rho_j^{\prime}$$

Vise versa,  for a surjective map $f: X \rightarrow C$ on a
curve with $rk H^1(C) \ge 2$ and a character $\rho$ we have $rk
H^1(C, \rho)=rk H^1(X,f^*\rho)$ (cf. \cite{1},).
Since on a hyperbolic curve $H^1(C,\rho) \ne 0$,  
the pullbacks of the local systems from $C$ belong to  
a component of $C^1_{\pi_1(X)}$ containing the identity
and in fact fills it. Moreover, the component $\rho_i f^*_i(H^1(C_i))$
belongs to $C^{-e(C)}_{\pi_1(X)}$ and has the dimension equal to 
$h^1(C)$.
Different, modulo automorphisms of the target,
admissible $f$ (and $C$) correspond to different components.
Hence the equivalences classes of admissible maps onto hyperbolic curves are in
one to one correspondence with components containing the identity
(here two maps are  equivalent, if they differ by an automorphism
 of the target). In particular one has enumeration of the 
holomorphic
maps from $X$ onto all possible hyperbolic curves in terms of
$\pi_1(X)/\pi_1(X)''$.

The group $Aut(X)$ of automorphisms of $X$ acts of $H^1(X,{\Bbb C}^*)=
Char(\pi_1(X)$ and hence on the components 
of $C^1_{\pi_1(X)}$. In particular, $\#Aut(X)=\nu_{C} \times \mu_{C}$ 
where $\nu_C$ is the number of components of $C^1_{\pi_1}$ containing 
the identity and belonging to the orbit corresponding to a map 
$X \rightarrow C$ and $\mu_C$ is the order of the stabilizer of this 
component in the action of $Aut(X)$ on components of $C^1_{\pi_1(X)}$.  

We shall conclude this section with two lemmas:

\vskip 0.3 cm






\proclaim{Lemma 1}
\roster
\item"{\bf a)}" $ \varphi_{\ast}:\pi_1(S) \to \pi_1(\Gamma)$ and
$ \varphi_{\ast}:h_1(S) \to h_1(\Gamma)$ are surjective;
\item"{\bf b)} " (\cite{24} Ch. 3, 11.4), 
\cite{25}, Ch. 4. Th.6): If $\ov S$ is of general type and $\ov\Gamma$
 is hyperbolic, then
   $$h_1(\ov F)\le\frac{e(\ov S)}{|e(\ov\G)|}+2.\tag 6$$
\endroster
\ep

\demo{Proof} a) is well known. For b) see (\cite {24} Ch. 3, 11.4), 
\cite{25}, Ch. 4. Th.6)

\hfill\quad\qed\enddemo

Note that the estimate (6) is sharp: if $\ov S$ is a
product of two hyperbolic curves, we have the equality.

\proclaim{Lemma 2}  
If the surface $S$ is affine, then 
$$|e(F)|\le\frac{ e(S)}{|e(\G)|}.$$\ep

\demo{ Proof of Lemma 2}

We may assume that $D$ ( in notations of  Lemma 1) intersects 
the  general fiber normally.
So, 
 outside the singular fibers  the map $\vp$ is
 locally 
trivial fiber bundle as well. 
Thus, we may use the standard procedure (\cite {24} (Ch. 3, 11.4), 
\cite{25}, Ch. 4. Th.6):

Let $ F$ be a general  and
 $ F_1\dots  F_s$ -- the singular
 fibers of the map $\vp$.
Then
$$e( S)=e( F)\cdot e(\G)+
\sum_{i=1}^s(e( F_i)-e( F)),\tag7$$

Similarly to projective case, since we assume that 
$F$ is connected, we have 
$$e( F_i)\ge e( F).\tag8$$ 
Inequality (8) for the affine surfaces was proven by M. Suzuki (\cite26)
and later, with more details, by  M. Zaidenberg (\cite{27}).

Since $ S$ is of general type, $ F$ should be hyperbolic
and  (8),(7) imply
$$e( F)\cdot e(\G)=e( S)-\sum_{i=1}^s(e(\tilde F_i)-e(\tilde F))\ge 0.$$

Thus, $$|e( F)|\le\frac{e( S)}{|e(\G)|},$$
and        $$h_1( F)=|e( F)|+1\leq \frac{e( S)}{|e(\G)|}+1.$$
\hfill\quad\qed\enddemo

\subheading{3. Projective case}

\vskip 0.3 cm

In this section we describe a bound for the number of maps
onto projective surfaces with $n(S)>0$ ($n(S)$ is the 
number of components of characteristic variety containing the
identity (cf. sect.1).
We shall put $T(h)=T_c(h/2)$ ( see (2)).

\proclaim {Definition} We call two maps $f_1:X\to Y_1$ and
$f_1:X\to Y_2$ equivalent if there exists an isomorphism
$\vp :Y_1\to Y_2$ such that the diagram
$$\alignat3
& &&\ X \\
&\qquad {}^{f_1}\swarrow && && \searrow{}^{f_2}\\
&\quad Y_1\qquad && \underset \vp\to\longrightarrow  &&\qquad Y_2
\endalignat$$
is commutative.\ep

\proclaim{Theorem 1}
Let $S$ be a smooth projective surface.
Then for the number $R(S)$ of possible rational dominant non-equivalent
maps from $S$
 onto all minimal surfaces $S_1$ of general
type with
$n(S_1)>0,$ the following inequality is valid:
$$R(S)\le \sum_{i=1}^{n(S)}T(a_i)T \left (\frac{e(S)}{a_i-1}+2\right )\le
n(S) T(h_1(S))T(e(S)+2),$$
where $a_i$  stands for dimension of $i-$th component of characteristic
variety  of $\pi_1(S).$
\ep

\remark{\bf Remark} If $n(S)>0$ for a smooth surface $S,$ then $n(S')>0$
for any other smooth surface $S',$ which is
birational to $S.$
\endremark
\proclaim {Lemma 3} Let $f_1:S\to S_1$ and  $f_2:S\to S_2$
be two rational  dominant maps from a smooth surface $S$
onto smooth projective minimal surfaces of general type $S_1,S_2$ respectively.
Let $\vp_1:S_1\to\G, \ \vp_2:S_2\to\G$ and $\vp:S\to\G$ be morphisms onto a
projective hyperbolic curve $\G$ such that the diagram
$$\alignat5
& S_1&&\overset{f_1}\to\longleftarrow &&\ S
&&\overset{f_2}\to\longrightarrow
&& S_2\\
&  &&{}^{\varphi_1}\searrow\ &&\ @VV \varphi V \ \swarrow{}^{\varphi_2}\\
 & && &&\ \Gamma \endalignat$$
is commutative. Let $G_{\g}=\vp_1^{-1}(\g), \
 H_{\g}=\vp_2^{-1}(\g), \ F_{\g}=\vp^{-1}(\g).$
 Assume that there is a Zariski open set $U\subset \G$ such that for each $\g\in U$
there exists an isomorphism $h_{\g}:G_{\g}\to H_{\g}$ and the diagram
$$\alignat3
& &&\ F_{\g} \\
&\qquad {}^{f_1}\swarrow && && \searrow{}^{f_2}\\
&\quad G_{\g} \qquad && \underset h_{\g}\to\longrightarrow  &&\qquad H_{\g}
\endalignat$$
is commutative as well. Then $S_1\cong S_2$ and the maps $f_1,f_2$
are equivalent.
\ep

\demo{Proof}
 Our aim is to obtain an isomorphism $g$ between $S_1$ and $S_2$
such that the diagram 
$$\alignat3
& &&\ S \\
&\qquad {}^{f_1}\swarrow && && \searrow{}^{f_2}\tag 9\\
&\quad S_1 \qquad && \underset g\to\longrightarrow  &&\qquad S_2,
\endalignat$$
is commutative.

Let $\gamma _1, \dots,   \gamma _n$ be the images in $\G$ of the 
indeterminacy points
 of the maps $f_1$ and  $f_2.$  Resolving the singularities of
$f_1$ and  $f_2$ we obtain a new surface which we also 
denote $S$. The new map $\vp$ takes each exceptional component 
into a point which is one of the points $\gamma _1, \dots,   \gamma _n$.
In particular the exceptional components would not
 intersect a general fiber.
It follows that we may
 assume that $f_1,f_2$ are morphisms (resolving
singularities of both of them).

Let $W$ be the image of $S$ in $S_1 \times_{\G} S_2$.
Consider the following diagram:
$$\alignat5
& &&\quad \quad &&\ S\\
&  &&{}^{f_1}\swarrow\ &&\ @VV \pi V  \ \searrow{}^{f_2} \\
& S_1&&\overset{ f'_1}\to\longleftarrow &&\ W &&\overset{
f'_2}\to\longrightarrow
&& S_2\\
&  &&{}^{\varphi_1}\searrow\ && && \ \swarrow{}^{\varphi_2}\\
 & && &&\ \Gamma \endalignat$$
 where 
$f_1', \ f_2'$ are induced rational maps. 
Since 
in the open set $\vp^{-1}(U)$ the maps $f_1'$ and $f_2'$
are one-to-one 
in $\pi(\vp^{-1}(U))$
(due to the  isomorphism on the general fiber), 
it follows that $S_1$ is birational to $S_2.$
But both are minimal, i.e. $S_1\cong S_2.$
 Moreover, for $g= f_2'\circ {(f_1')}^{-1}$ diagram (9) is commutative.
\hfill\quad\qed\enddemo

\demo{Proof of Theorem 1}
Let $S_i, \ i=1, \dots, N$ be  surfaces of general type with $n(S_i)>0.$
Let $\vp_i:S_i\to\G_i$ be regular maps with a connected general fiber
 onto a hyperbolic curve and $f_i: S\to S_i$ be a rational
 dominant map.
Consider the commutative diagram:
$$\alignat3
& &&\ \tilde S \\
&\qquad {}^\pi\swarrow && && \searrow{}^{\tilde f_i}\\
&\quad S\qquad && \underset f_i\to\longrightarrow  &&\qquad S_i\tag10\\
&{}^{\sigma _i}\downarrow && &&\qquad \downarrow{}^{\vp_i}\\
&\quad \Sigma _i\qquad && \underset {\vp_{f_i}}\to\longrightarrow
 &&\qquad \G _i
\endalignat$$
where $\pi: \tilde S\to S$ is a resolution of singularities of all
the maps
 $f_i$ (i.e., $\tilde f_i$ are regular on $\tilde S),$
$\Sigma _i$ stands for   the Stein factorization of the map
 $\vp_i\circ \tilde f_i.$
Since $\G_i$ is hyperbolic, for any exceptional curve $E$ of
the map $\pi$ the dimension  
$\dim\vp_i\circ \tilde f_i(E)$ is zero, i.e., the map $\sigma _i$ is
everywhere defined on $S.$

If two maps $f_1,f_2$ are not equivalent, then
(due to Lemma 3)

-- either the maps $\vp_{f_i}\circ \sigma_i,  \ i=1,2$ are not equivalent;

--or  the maps $\vp_{f_i}\circ \sigma_i,  \ i=1,2$ are  equivalent, but 
the restrictions of $f_1$ and $f_2$ on a general fiber $F_{\sigma}=
\sigma_i^{-1}(s), \ i=1,2$
are not equivalent for a general $s\in \Sigma_1$
(in this case $\Sigma_1=\Sigma_2,$ since they are the Stein 
factorization of the same map).

There may be at most $n(S)$ (up to isomorphism)
 curves $\Sigma _i,$ which may appear in
diagram (10), at
 most $T(h_1(\Sigma _i))$ non-equivalent maps $\vp_{f_i}$ for each
$\Sigma _i,$ and
at most $ T(h_1(F_{\sigma_i}))$ non-equivalent maps of the fiber.
Thus,
$$ R(S)\leq\sum_{\Sigma _i}T(h_1(\Sigma _i))T(h_1(F_{\sigma_i})).$$

Here $h_1(\Sigma _i)=a_i$  is precisely
the dimension of the corresponding component of
characteristic variety of $\pi_1(S)$
and estimate for $ h_1(F_{\sigma_i})$
is provided by Lemma 1.

Thus $$R(S)\le \sum_{i=1}^{n(S)}T(a_i)T \left
 (\frac{e(S)}{a_i-1}+2\right ).$$

Due to the same Lemma 1 we have $2\le a_i\le h_1(S)$, i.e.
$$\sum_{i=1}^{n(S)}T(a_i)T \left
 (\frac{e(S)}{a_i-1}+2\right )\le
n(S) T(h_1(S))T(e(S)+2),$$
 \hfill\quad\qed\enddemo

\vskip 0.3 cm
\proclaim {Remark} An alternative way to estimate the
number of pairs $(\Sigma_i,\G_i)$ is the following.
  Since the surface $S$ is of general type, $e(S)\ge e(S_m),$
where $S_m$ is its minimal model.  The linear system $|5K(S_m)|$
is base point free and provides the birational map (\cite{28}).
Hence the number of non-equivalent maps $S\to \G_i$
does not exceed the same number for the general divisor
$L\in |5K(S_m)|.$
$$g(L)=\frac {(L,L+K(S_m))}{2}+1\le 15(K(S_m))^2+1\le 45e(S_m)+1\le45e(S)+1.$$
Therefore,
$$R(S)\le T(45e(S)+1) T(e(S)+2).$$

This estimate is obviously worse then one, obtained in Theorem 1.
For example, if $n(S)=0,$ then $R(S)=0.$
Moreover, $h_1(S)\le h_1(L)\le 45e(S)+1.$
\ep

\vskip 0.3 cm
\proclaim{Theorem 2} Let $X\subset \Bbb P^N$ be a smooth projective variety
with nef and big canonical bundle and $n=dim X.$ Let $s\ge 1$ be a rational
 number such that

1) $sK(X)$ is a Cartier divisor on $X$;

2) linear system $|sK(X)|$ is base point free.

Then $$ R(X)\leq n(X)T(h_1(X))T(5s^nK(X)^n+38),$$
where $R(X)$ is, as above, the number of dominant maps from $X$ onto
minimal surfaces $S$ of general type with $n(S)>0.$ 
\ep

\vskip 0.3 cm
\proclaim{Remark} Due to Theorems of V.Shokurov and Y. Kawamata
(\cite{29}, \cite{30}) such $s$ always exists. J.Kollar
(\cite{31})
proved that any integral  $s>2(n+2)!(2+n)$ will do.
If $K(X)$ is ample, then any $s>1/2(n(n+1))$ will be good (\cite{32}).
If $dim X=3$ then $s=7$ (\cite{33}). Theorem 2 may be proved
without much changes for varieties, having at most  canonical
singularities ( then the estimate would depend on the index of $X$),
 but then we do not know the good effective estimate of $s.$\ep

  \demo {Proof of Theorem 2} We use induction on $n=dim X.$

Let $dim X=2.$ Then, by Theorem 1

$$R(S)\le n(S) T(h_1(S))T(e(S)+2) .$$
 Since for minimal surfaces of general type
$$e(S)\leq 5K(X)^2+36$$
(\cite{24}, Corollary 3.2, ch. VII),
we have:
$$R(S)\le n(S) T(h_1(S))T( 5K(S)^2+38) ,$$
and the claim follows for $s\geq 1.$

\vskip 0.3cm

Now let $dim X=n>2.$ By Bertini Theorem, the general
divisor $L\in|sK(X)|$ is a smooth
subvariety of $X.$ It has the following properties:

1) $dim L=n-1;$

2) $K(L)=(K(X)+L)\bigm|_L=(s+1)K(X)\bigm|_L;$

3) $K(L)$ is nef and big;

4) $[K(L)]^{n-1}= ((K(X)+L)^{n-1},L)=(s+1)^{n-1}s[K(X)]^{n};$

5) $\frac{s}{s+1}K(L)=[\frac{s}{s+1}(s+1)K(X)]\bigm|_L=sK(X) \bigm|_L$
  is a Cartier divisor;

6) $\bigm |\frac{s}{s+1}K(L)\bigm |$
    is a base point free linear system.

 Let $\Vp : X\to \Bbb P^M$ be  morphism,  defined by the
 linear system $|sK(X)|$
and $\ov\Vp :\Bbb P^N\to\Bbb P^M$ some its rational extension into
 projective space.  By construction, $L= \Vp^{-1}(H\cap\Vp(X))\cap X=
 H'\cap X, $ where $H$ is a hyperplane in $\Bbb P^M$ and
 $H'=\ov\Vp^{-1}(H),$ i.e.
 $H'$ is a
hypersurface (smooth for a general $H$ due to Bertini Theorem).

 Since a smooth hypersurface  in $ \Bbb P^N$ is its hyperplane
section for some Veronese embedding, the Lefschetz Theorem on
 hyperplane section (\cite{34})
 is valid for it , and
$$\pi_1(L)=\pi_1(H'\cap X)\cong\pi_1(X).$$

Hence

   7) $n(L)=n(X).$

The inductive assumption yields:
 $$ R(L)\leq n(L)T(h_1(L))T\big (5{\big
 (\frac{s}{s+1}\big )}^{n-1}K(L)^{n-1}
+38 \big )=$$
  $$ n(X)T(h_1(X))T(5s^nK(X)^n+38).$$

Since $R(L)\geq R(X)$ for a general $L,$ the Theorem is proved.
\hfill\quad\qed\enddemo

\subheading{4. Maps into affine curves and  surfaces of general (hyperbolic)
 type}

\vskip 0.3 cm
 Let  $T'(h)=max\{T(g,r), \ \text{with} \quad 2g+r-1=h\}$ (see(3)).
\vskip 0.3 cm

\proclaim{Theorem 3}
Let $X$ be a quasiprojective variety and $N(X)$ -- the number of all regular
dominant maps from
$X$ into hyperbolic curves.
Then
$$N(X)\le \sum_{i=1}^{n(X)} T'(a_i)\le  n(X)\cdot T'(h_1(X)),$$
where $a_i$ is the dimension of the $i-$th component of the
characteristic variety of $\pi_1(X).$
\ep

\demo{Proof of Theorem 3}
Let $X\overset \varphi\to\rightarrow \Gamma$ be a dominant regular map
 from $X$ into a
hyperbolic curve $\Gamma.$
Let $\ov X,\ov\G$ be closures of $X$ and $\G$, respectively,
 such that the map $\vp$ may be
extended to a regular map $\ov\vp: \ov X\to \ov\G.$

Consider the Stein factorization
of $\ov\vp,$ i.e.   a smooth projective curve $\tilde\G$
and   a map $\tilde\vp: \ov X\to\tilde \G,$ such that in the commutative
 diagram
$$\alignat 3
& \ov X\quad&&\overset {\ov\vp}\to\longrightarrow &&\quad\ov\G\\
&{}^{\tilde\vp} \searrow &&  &&\nearrow {}^h
\quad\\
& &&\quad\tilde \G
\endalignat$$
the regular map $\tilde\vp$ has connected fibers, and $h$ is finite.
Let $\tilde \G'=\tilde\vp(X).$ We have
$$N(X)\le\sum_{\tilde \G'}T'(h_1(\tilde \G'))\le\sum_{i=1}^{n(X)} T'(a_i).$$

Since the fibers of $\tilde \vp $ are connected, $h^1(\tilde\G)\le h^1(X)$
(see Lemma  1), and
 the number of possible $\tilde\G'$  is equal to $n(X),$
 we have
$$N(X)\le n(X)\cdot T'(h^1(X)).$$\hfill\quad\qed
\enddemo

\proclaim{Remark} Both numbers depend only on the fundamental group of $X.$
\ep
\proclaim{Remark}
In particular, this Theorem provides  positive answer
 to an old conjecture of V.Lin and E.Gorin (\cite{22},4.2).

{\bf Conjecture}(E.Gorin, V.Lin) The number $\sigma(X)$
of all holomorphic functions omitting two values on affine
 variety $X$ may be bounded from above by a function, depending on
topology of $X$  only.

This conjecture is true for curves. In \cite{20} a bound 
was found for $\sigma(X),$ depending on $h_1(X)$ and $h_2(X),$
when $dim X=2.$
Theorem 3 provides the following estimate.\ep

\proclaim {Corollary 1}
$$\sigma(X)\leq \sum_{i=1}^{n(X)} (\sum_{2g+r-1=a_i}T_0(g,r))\le
N(X)$$
has a bound, depending
 on the fundamental group of $X.$
\ep

Now we shall estimate the
number of maps between two affine surfaces of general type
with $n(S)>0$ and meeting some additional conditions.
In particular, we are going to estimate
$\#Aut (S)$ of all the
 automorphisms of this class of affine surface $S.$
We want to use the fact that if $n(S)>0,$ then the surface is fibered.

The idea is the following one.
Consider all the maps, preserving a fibering.
If the map preserve fibering, it maps a
fiber onto itself and provides
an automorphism of the fiber.
If two maps coincide on the open set, they coincide.
That means that two different maps may coincide only
on the finite number of fibers.
Since the number of maps is finite,
in the general fiber all the maps are different.
Hence, the number of such automorphisms do not
 exceed he number of automorphisms of a general fiber.
Further on we denote by $L(C)$ the number of automorphisms
of a curve $C.$

\proclaim{Theorem 4}
   For smooth affine surface $S$ of general type with $n(S)\ge 0,$
 $\#Aut(S)$ has an estimate
$$\#Aut(S)\le\min_{\nu _{\G}\ne 0}\{\nu _{\G}\cdot
L(\Gamma)\cdot L_{\vp}\}.\tag11$$
where   $\nu _{\G}$ is the length of an orbit of a component 
of characteristic variety corresponding to a map $\vp: 
S \rightarrow {\G}$, $F_{\vp}$ is a generic fiber of  $\vp$ 
and $L_{\vp}= \max L(F_\vp)$ for all such $\vp.$

This estimate
may be simplified as
$$\#Aut(S)\le \min_{n_i(S)\ne 0} \{252 n_i(S)\cdot e(S)\}.\tag12$$
\ep

\demo{Proof of Theorem 4}
Let $\vp:S\to\Gamma$ be a dominant map of $S$ (with a connected general fiber)
into a hyperbolic curve $\G,$ \
$f: S\to S$  an automorphism.

We have the following diagram:
$$\alignat3
&\quad S\           && \overset{ f }\to\longrightarrow\ &&\ \ S\\
&{}^{\vp}\downarrow &&                           &&{}^{\vp}\downarrow\tag13\\
&\quad \G\ &&\overset{\vp_ f}\to{\cdots>}\ &&\ \ \Gamma
\endalignat$$

Since $f$ is an isomorphism, the general
fibers of $\vp$ and $\vp\circ f$ are
 both connected. That's why if $\vp_f$ is defined, it has to be an isomorphism,
i.e. $\vp$ and $\vp\circ f$ define the same irreducible component
of characteristic variety of $\pi_1(S).$ Then $\vp $ induces an isomorphism
of $F_\vp=\vp^{-1}(\g)$ onto $\vp^{-1}(\vp_f(\g))$ for
a general point $\g\in \G.$ Thus, there may be at most
$$ a_{\Gamma}=L(\G)\cdot L(F_\vp)$$ automorphisms $f$ for which $\vp_f$
is defined.

Assume that diagram (13) cannot be completed by the map $\vp_f,$ that is,
$\vp$ and $\vp\circ f$ define different
components (of the same dimension) of the characteristic variety
of $\pi_1(S).$
Let $f_1,f_2$ be two  automorphisms such that $\vp\circ f_1,$ and
$\vp\circ f_2$ define the same component.
Then changing $f$ for $f_2\circ f_1^{-1}$ in
diagram (13), we obtain the new diagram (13), which may be completed by
the isomorphism of the curve $\G.$

Therefore,  $$\#Aut(S)\le\nu _{\G}\cdot L(\Gamma)\cdot  L_{\vp},$$
We remind that here    $\nu _{\G}$ stands for the number of
 the components of characteristic
variety of $S,$ which correspond to (non-equivalent) maps $\vp:S\to \G$
belonging to the same orbit, 
and $L_\vp= \max L(F_\vp)$ for all such $\vp.$

Since this is valid for any $\G,$ provided that a dominant map $S\to\G$
exists, we have
$$\#Aut(S)\le\min_{\nu _{\G}\ne 0}\{\nu _{\G}\cdot
L(\Gamma)\cdot L_{\vp}\}.$$

 $\G$ may be  $a priori$ compact, that is why the best estimate is
 $L(\G)\le 42|e(\G)|$. On the other hand, $F_\vp$ is an affine curve
and (due to Lemma 2)
 $$|e(F_\vp)|\le\frac{ e(S)}{|e(\G)|}.$$

Hence, (see(1)),
 $$L(F_\vp)\le 6\frac{ e(S)}{|e(\G)|}.$$

Inserting these estimates into (11), we get (12).
\hfill\quad\qed\enddemo

\proclaim{Remark} The second bound in this theorem is not sharp,
 but the bound (11) is. In the Appendix we 
 demonstrate this by  examples.\ep

\proclaim {Theorem 5}  Let $S$ be an affine surface
of general type.
Let $S'$ be any other surface of general type with $n(S)>0.$
Then the number $R(S,S')$ of dominant morphisms
from $S$ into $S'$
$$R(S,S')\le n(S)T'(e(S)+1).$$
\ep

\demo{Proof of Theorem 5}
Let $f:S\to S'$ be a dominant morphism, $\vp'$ a map with connected
 general fiber from $S'$ onto a hyperbolic curve $\G'.$
Denote:

-- $\ov S'$ such a closure of $S'$ that the extension $\ov \vp'$ of $\vp'$
onto $\ov S' $ is a morphism;

-- $\pi:\ov S\to \tilde S$ such a resolution of $\tilde S$, that the extension
 $\ov f$ of $f$ onto $\ov S$ is a morphism;

-- $\ov\G'$ - the closure of $\G';$

-- $\ov\G$ -the Stein factorization of $\ov\vp'\circ\ov f;$

-- $\G=\vp(S)\subset\ov\G.$

The following  diagram is commutative:
$$\alignat6
& \tilde S  && \ @<<\pi< &&\overline S &&@>\ov f>>  &&\overline  S'\\
&\cup   && &&  @VV\overline\varphi  V&&     @VVV \overline\varphi'\\
&S && && \overline \G && @>>>  &&\overline \G'\\
  & &&\varphi\!\!\!\searrow   \quad && \cup && &&\cup\\
& && &&\Gamma &&@>>\varphi_f>  &&\Gamma'
\endalignat$$

There may be at most $ n(S)$ ways to
 choose the map $\vp$ and ( $\ov\vp$) in this digram.
For each choice of $\vp$ may be at most $T'(h_1(F))$ maps
from a fiber $F=\vp^{-1}(\g)$ over a general point $\g\in\G$
into a fiber $(\vp')^{-1}(\vp_f(\g)).$
As it was shown in Corollary 1
  $$h_1( F)=|e( F)|+1\leq \frac{e( S)}{|e(\G)|}+1.$$

Hence, $$R(S,S')\le n(S)T'(e(S)+1)$$
\hfill\quad\qed\enddemo

\subheading{{\bf Appendix: \ Automorphisms of arrangements}}

\bigskip

Let $M=\cupl_{i=1}^nL_i$ be an arrangement of $n$ lines in
$\Bbb P^2$, $S=\Bbb P^2-M$ and
$A_1,\dots,A_t$ be the intersection points of $d(A_i)\ge3$ lines from $M.$
 If $t>0$, then projection from $A_1$ 
yields  a regular map
$\varphi_1: S\to\Bbb P^1-\{z_1,\dots z_d\}$, where
$d=d(A_1),$ of $S$ onto a hyperbolic curve.
If $S$ is of general ( hyperbolic) type, the fibers should
 be hyperbolic as well.
Thus, we may
use the above results to estimate the number of automorphisms of $S.$
 According to  (1) and (11)
 (we take minimum here over a smaller set, than in
(11))
$$\#Aut(S)\le\min_{0\le k\le t} \{n_{d_k-1}(S)A(0,d_k)A(0,n-d_k+1)\}\le$$
$$\min_{0\le k\le t} \{25n_{d_k-1}(S)d_k(n-d_k+1)\}\le
\frac{25}{4}n(S)(n+1)^2.$$

{\bf Example 1.} Consider the arrangement of 6 lines, which are
 the sides and medians of a triangle.

Let $\{w_1:w_2:w_3\}$ be coordinates in $\Bbb P^2,$
$L_1=\{w_1=0\},$
$L_2=\{w_2=0\},$
$L_3=\{w_1=w_2\},$
$L_4=\{w_1+w_2-w_3=0\},$ $L_5=\{w_1+2w_2-w_3=0\},$ $L_6=\{2w_1+w_2-w_3=0\},$
$S=\Bbb P^2-\cupl_{i=1}^4L_i.$

The vertices of the triangle: $A_1=(0:0:1), \ A_2=(1:0:1), \  A_3=(0:1:1),$
and the intersection point of medians:
$ A_4=(1:1:3)$
have $d_i=3,  \ i=1,..., 4.$  Linear projection from each $A_i$ 
yields the map $\vp_i$ of $S$ onto 
$\G=\Bbb P^1-\{0,1,\infty\}$.
There is also the fifth map onto
 $\G:$
$$\vp_5=\frac{w_1(2w_1+w_2-w_3)}{w_2(w_1+2w_2-w_3)}.$$

A calculation with the fundamental group as in section 2 
or the results in \cite{3}
yield 
$n(S)=n_2(S)=5.$

The first component is defined by the fibering $\vp_1:S \to  \G,$
where  $\vp_1=(w_1:w_2).$
The general fiber, being $\Bbb P^1$
with four general punctures has
four automorphisms, i.e.  $L_{\vp_1}=4,  \ L(\Gamma)=A(0,3)=6. $
Let us show, that all these  automorphisms may be realized as the
 automorphisms of the surface
$S.$ 

Let  $k=w_1/w_2,$ \ $t=(w_3-w_1-w_2)/w_2.$ In these coordinates
 $S= \{(k,t)\in \Bbb  C^2:(k,t)\not\in \cup L'_i\},$ where
$L'_1=\{k=0\},$  \ $L'_2=\{k=1\},$
$L'_3=\{t=0\},$ \ $L'_4=\{t=1\}, $ \ $L'_5=\{t=k\};$
 \ $\vp_1(k,t)=k;$ \
$\vp_1^{-1}(k)=F_{\vp_1}=\{t\ne 0,1,k\}. $
Four automorphisms of the fiber are induced by the following
 automorphisms of $S:$
$$f_1(k,t)=(k,t),  \ \  f_2(k,t)=(k,k/t),
\  f_3(k,t)=(k,\frac{k-t}{1-t}), \  \  f_4(k,t)=(k,\frac{k(t-1)}{t-k}).$$
Six  automorphisms of the base  $\G=\Bbb P^1-\{0,1,\infty\}$ 
may be realized as
 $g_1(k,t)= (1/k,1/t),$
 $g_2(k,t)= (1-k,1-t),$
 $g_3(k,t)= (1-1/k,1-1/t),$
$g_4(k,t)= (\frac{1}{1-1/k},\frac{1}{1-1/t}),$
$g_5(k,t)= (\frac{1}{1-k},\frac{1}{1-t}),$
$g_6(k,t)= (k,t).$

It follows that  we have  
 $L(\G)\cdot L_{\vp_1}=A(0,3)\times 4=6\times 4=24 $ 
automorphisms of the surface $S,$ corresponding to the map $\vp_1.$

  Now we want to show that all five components belong to the same orbit.
In other words, we want to 
 ``connect'' the different components of
 the characteristic variety by automorphisms, i.e.
to find $f_i \in Aut(S),$ such that
 $\vp_1\circ f_i$ is $\vp _i, \ i=2, ..., 5$ ( see digram (13)).

 The map $\a :\Bbb P^2
\to  \Bbb P^2, \ \a(w_1:w_2:w_3)=(w_2: -(w_1+w_2-w_3):w_3),$
 permutes the vertices of the triangle, and leaves the
 point $A_4=(1:1:3)$  fixed. Hence, it preserves the arrangement
and ``connects'' the components,  defined by $\vp_1,\vp_2,\vp_3.$

The map $\be :\Bbb P^2
\to  \Bbb P^2, \ \be (w_1:w_2:w_3)=(2w_1+w_2-w_3:w_1-w_2:-w_3+4w_1)$
 sends $A_4$ to $A_1,$ preserving the arrangement, i.e. ``connects'' the
fourth and the first components.

The map $\g:\Bbb P^2
\to  \Bbb P^2,$  \ 
$\g(w_1:w_2:w_3)=(w_1(w_1+2w_2-w_3):w_2(2w_1+w_2-w_3):
w_1(w_1+2w_2-w_3)+w_2(2w_1+w_2-w_3)+w_1w_2)$
connects the first and the fifth component:
$$\g(L_1\cup L_5)=L_1, \ \g(L_2\cup L_6)=L_2, \ \g(L_3\cup L_4)=L_3,
 \ \g(A_1)=L_4,  \ \g(A_2)=L_5, \ \g(A_3)=L_6.$$
$ \g(L_1')=A_3, \ \g(L_2')=A_2, \ \g(L_3')=L_3, \  \g(L_4')=A_4,
\ \g(L_5')=L_1, \ \g(L_6')=L_2, \ \g(A_1)=L_4, \ \g(A_2)=L_5,  \ \g(A_3)=L_6.$

Altogether,  $\# Aut(S)=n_2(S)A(0,3)L_{\vp_1}=5\times 6\times4 =120,$
which shows that the estimate (11) is sharp.

\bigskip
{\bf Example 2.} Consider the following arrangement of 
 $12$ lines and $9$ points:
$$L_1=\{w_1=0\}, \
L_2=\{w_2=0\}, \ 
L_3=\{w_3=0\},$$
$$L_4=\{w_1+w_2+w_3=0\}, \ 
L_5=\{w_1+ew_2+e^2w_3=0\}, \ 
L_6=\{w_1+e^2w_2+ew_3=0\},$$
$$L_7=\{w_1+w_2+ew_3=0\}, \ 
L_8=\{w_1+ew_2+w_3=0\}, \ 
L_9=\{w_1+e^2w_2+e^2w_3=0\},$$
$$L_{10}=\{w_1+w_2+e^2w_3=0\}, \ 
L_{11}=\{w_1+ew_2+ew_3=0\}, \ 
L_{12}=\{w_1+e^2w_2+w_3=0\}.$$

 Let $N=\{(1:1), (0:1), (e:1), (e^2:1)\},$ where
$ e^3=1$ be the set of 4 points in $\Bbb P^1$.
 There is one component (see\cite{3}), corresponding to the map $\vp$
of all
the arrangement to $\Bbb P^1- N.$
This map sends the point on the cubic
$$ a (w_1^3+w_2^3+w_3^3)+bw_1w_2w_3=0$$ to $(-3a,b).$
The fibers are  general tori with 9 punctures, so
$$\nu_{\Gamma}=1,  \ L(\Gamma)=12,  \  L_{\vp}=18, \ 12\times 18=216$$
 Thus, in this case formula (11) gives  at most 216, which is sharp
(\cite{35} p. 298)

\bigskip
{\bf Example 3.} Consider the arrangement dual to the one in 
example 2:  9 lines  $L_i, i=1,...,9$
intersecting in 12 points $A_j, j=1,...,12.$  At each point
three lines are intersecting.
$$L_1=\{w_2-w_3=0\},
 \  L_2=\{w_1-w_3=0\},
 \  L_3=\{w_1-w_2=0\},$$
 $$L_4=\{ew_2-w_3=0\},
 \ L_5=\{ew_1-w_3=0\},
  \ L_6=\{ew_1-w_2=0\}$$
 $$L_7=\{e^2w_2-w_3=0\},
 \ L_8=\{e^2w_1-w_3=0\},
 \ L_9=\{e^2w_1-w_2=0\}.$$
$$  A_1=(1:0:0), \   A_2=(0:1:0), \   A_3=(0: 0:1), \ A_4=(1:1:1),$$
 $$  A_5=(1:1:e), \ A_6=(1:1:e^2),
  A_7=(1:e:e^2), \   A_8=(1:e:1), $$  $$ A_9=(1:e:e), \ A_{10}=(1:e^2:e),
 \  A_{11}=(1:e^2:e^2), \ A_{12}=(1:e^2:1).$$

Let $S=\Bbb P^2- \cup L_i.$
There are $16$ $2$-dimensional components of characteristic 
variety of $S$ ( \cite{3}): 12 maps $\vp_i$ to $\G=\Bbb P^1-\{0,1,\infty\},$
with the only singular point at each $A_i, i=1,...,12$ and 
four  maps $\psi_i, i=1,...,4$, to $\G'\sim \G,$ 
 each  for  a choice of 3 points
which are not connected by lines. The fibers of $\vp_i$ are lines, 
the fibers of $\psi_i$ are elliptic curves, thus  the components 
corresponding to $\vp_i$ and  $\psi_i$ belong to different orbits.
 On the other hand, one can write explicitly  12 linear automorphisms
of $S,$ permuting the points $A_i.$  Hence,  $\nu_{\G}= 12.$

Consider a fiber $ F_a=\{\vp_1=a\}=\{(w_1:w_2:w_3): w_2-w_3=a(ew_2-w_3)\}.$
It is isomorphic to $C= \Bbb P^1-\{0,1,t, e,et, e^2, e^2t\},$
where $e^3=1$ and $t=\frac{ea-1}{a-1},$
and $L(C)=3.$
Thus  the final estimate is $12\times 6\times 3=216,$
which is an expected estimate, because this configuration of points
and lines is dual to  one from example 2.

\subheading{Acknowledgments}
The  authors are   grateful to S. Yuzvinsky for his help   
and to V. Lin 
for providing the most  useful information and references.

\end
\Refs
\widestnumber\key{34}

\ref\key{1}\by D. Arapura\paper Geometry of cohomology
support of local systems I,\jour J.Algebraic Geom.
\vol 6\yr 1997\pages
563-597\endref


\ref\key {2}\by A.Beauville \paper Annulation du $H^1$ pour les fibr\' {e}s
en droites plats\inbook Lecture Notes in Math.\vol 1507 \publ
Springer-Verlag\yr 1992\pages 1-15\endref

\ref\key{3} \by A.Libgober \paper Characteristic varieties of algebraic
curves, \jour Applications of algebraic geometry to coding theory, physics
and computations (Eilat, 2001), NATO Sci. Ser. II, Math. Phys. 
Chem., 36, Kluwer Acad. Publ., Dordrect, 2001. \pages 215-254 \endref

\ref \key{4}\by S. Kobayashi, T.Ochiai \paper Meromorphic
mappings into complex spaces of general type \jour Inv. Math. \vol 31
\pages 7-16 \yr 1975 \endref


\ref\key{5}\by S. Iitaka\paper On logarithmic Kodaira dimension of
algebraic varieties \inbook Complex Analysis and Algebraic Geometry\publaddr
Tokyo \publ Iwanami \yr 1977 \pages 178-189\endref



\ref\key{6}\by F.Sakai\paper Kodaira dimension of complements of divisors
\inbook Complex Analysis and Algebraic Geometry\publ
Cambridge university Press\yr 1977\pages 239-259\endref

 \ref\key{7}  \by T. Tsushima \paper Rational maps of varieties of
hyperbolic type \jour
       Proc. Jap. Acad. \vol ser A, 55 \pages 95-100 \yr 1979\endref


\ref\key{8}\by  K. Maehara \paper  A finiteness property of variety of
general type\jour  Math. Ann.\vol 262\pages 101-123 \yr 1983\endref



\ref \key {9} \by E. Szabo \paper Bounding the automorphisms groups
\jour Math. Ann.\vol 304 \yr 1996 \pages 801-811\endref


\ref\key{10}\by G. Xiao \paper Bound of automorphisms of surfaces of general
type,I \jour Ann. of Math.\vol 139\yr 1994\pages 51-77\endref

\ref\key{11}\by  G. Xiao\paper Bound of automorphisms
 of surfaces of general type. II. \jour J. Algebraic Geom \vol 4 \yr 1995
 \pages 701--793 \endref



\ref\key {12}\by R.D.M.Accola \paper
Topics in the theory of Riemann surfaces
\inbook Lecture Notes in Math.\vol 1595 \publ
Springer-Verlag\yr 1994\endref

\ref\key{13} \by A. Howard and A. Sommese \paper On the theorem of de Franchis \jour
Annali Scuola Norm. Sup. Pisa \vol 10 \yr 1983 \pages 429-436 \endref



\ref\key{14}\by  E. Kani \paper Bounds on the number of non-rational subfields
 of a function field \jour Inv. Math. \vol 85\pages 199-215 \yr 1986\endref




\ref\key {15}\by A.Alzati, G.P.Pirola
\paper Some remarks on the de Franchis Theorem
\jour Ann. Univ. Ferrara, Sez VII, Sc.Mat. \vol  36\yr 1990
\pages 45-52 \endref








\ref\key{16}\by T. Bandman, D. Markushevich \paper On the number of
 rational
 maps between varieties of general type\jour J. Math. Sci. Univ. Tokyo
 \vol 1 \pages 423-433 \yr 1994 \endref




    \ref\key{17}    \manyby     T.
Bandman \paper     Surjective holomorphic mappings   of
projective
    manifolds\jour  Siberian Math.      Journ.  \vol    22  \pages
204 -210 \yr 1982     \endref


 \ref\key{18} \by T. Bandman, G. Dethloff \paper Estimates of the
 number of rational
mappings from a fixed variety to  varieties of
general type \jour
Ann. Inst. Fourier \vol 47
\yr 1997 \pages 801-824\endref







\ref\key{19} \by T. Bandman \paper Topological invariants of a variety
 and the number of its
holomorphic mappings \inbook
 J. Noguchi (Ed.): Proceedings of the International
Symposium Holomorphic Mappings, Diophantine Geometry and Related Topics
\pages 188-202 \publ
RIMS \publaddr  Kyoto University \yr 1992 \endref








\ref\key{20} \by T. Bandman \paper Holomorphic functions minus two values
 on an affine surface \jour Vestnik Moskovskogo Univ. Matematika \vol35
\yr 1980 \pages 43-45 \endref


\ref\key{21} \by G. Heier \paper Effective finiteness  theorems
for maps between canonically polarized compact complex manifolds,
preprint, AG-0311086 \endref

\ref\key{22} \by M.G. Zaidenberg and V.Ya. Lin \paper Finiteness
theorems for holomorphic maps \inbook  Several Complex Variables III, Encyclopedia Math.
Sciences \vol  9 \pages 113-172 \publ  Springer Verlag \yr 1989 \endref


\ref\key{23} \by A.Libgober \paper Hodge decomposition of Alexander 
invariants. \jour Manuscripta Math. \vol 107 no. 2 \yr 2002 \pages 251-269
\endref

\ref\key{24} \by W. Barth, C. Peters, A. van de Ven \book Compact
complex surfaces \publ Springer Verlag  \yr 1984 \endref


 \ref\key  {25}
\by I. R.  \v Safarevi\v c    \book Algebraic Surfaces \publ (Proc.   of
the   Steclov                 Inst. of Math., {\bf 75} 1965), AMS,
Providence  \yr 1967 \endref


\ref \key {26} \by M. Suzuki \paper Sur les op\'erationes
holomorphes du group additif complexe
sur l'espace de deux variables complexes
\jour  Ann. Sci Ec. Norm. Sup., serie
4,\vol 10 \yr 1977 \pages 517-546\endref




\ref\key{27} \by M.G. Zaidenberg \paper "Isotrivial families of
curves on affine surfaces and the characterization of the affine plane,"
 \jour Math. USSR-IZVESTIYA \vol  30 \yr 1988 \pages 503-532 \endref



\ref\key{28}\by I. Reider\paper Vector bundles of rank 2
and linear systems on Algebraic surfaces \jour Ann. Math \vol 127 \yr
1988 \pages 309-316 \endref




\ref\key{29}\by V.Shokurov\paper The non-vanishing theorem\jour
Math.USSR-Izv.\vol 19 \yr 1985\pages591-607\endref


\ref\key{30}\by Y.Kawamata\paper On the finiteness of generators
 of a pluricanonical ring for a threefold of general type
\jour Amer.Jour. Math.\vol 106\yr 1984\pages 1503-1512
\endref














\ref\key{31} \by J. Kollar \paper Effective base point freeness
\jour Math.Ann.\vol 296 \yr 1993 \pages 595-605\endref


\ref\key{32}\by U.Angehrn, Y.-T.Siu \paper Effective freeness
and point separation for adjoint bundles\jour Invent. Math.\vol 122
\yr1995 \pages 291-308\endref


\ref\key{33}\by L.Ein, R.Lazarsfeld\paper Global generation of
pluricanonical and adjoint linear series on smooth projective
threefolds \jour J.Amer.Math.Soc.\vol 6\yr 1993\pages 875-903\endref


\ref\key{34}\by M.Goresky, R. MacPherson\paper Stratified Morse Theory
\publ Springer-Verlag \yr 1988\endref


\ref\key{35}\by E. Brieskorn, H. Knorrer \book Plane Algebraic curves
 \publ Birkhauser \yr 1986 \endref






























\endRefs
